\documentclass[12pt]{article}

\usepackage{graphicx} 
\usepackage{amsmath} 
\usepackage{amssymb,amsthm}
\usepackage{color}
\usepackage{tikz}
\usepackage{pgfplots}

\newcommand{\R}{{\mathbb R}}\newcommand{\N}{{\mathbb N}}

\newcommand{\Z}{{\mathbb Z}}\newcommand{\C}{{\mathbb C}}

\newcommand{\hZ}{\widehat{Z}}

\let\epsilon\varepsilon
\let\theta\vartheta

\newtheorem{theorem}{Theorem}[section]\newtheorem{lemma}[theorem]{Lemma}

\newtheorem{corollary}[theorem]{Corollary}

\newtheorem{remark}[theorem]{Remark}

\title{Validity of the stochastic Ginzburg-Landau approximation in higher space dimensions \\
-- A Wiener algebra approach --}
\author{Anna Logioti, Guido Schneider \\
{\small
Institut f\"ur Analysis, Dynamik und Modellierung, Universit\"at Stuttgart,  } \\ {\small  Pfaffenwaldring 57,  70569 Stuttgart, Germany}}

\begin{document}

\maketitle

\begin{abstract}
We consider an anisotropic $d$-dimensional Swift-Hohenberg model $ \mathcal{O}(\varepsilon^2) $-close to the first instability, where $ 0 < \varepsilon \ll 1 $
is a small perturbation parameter.
This model for pattern formation is perturbed 
with additive noise in time and space. 
By a multiple scaling ansatz we 
derive a stochastic $ d $-dimensional Ginzburg-Landau equation for the approximate description
of the bifurcating solutions.
We prove the validity
of the approximation by this amplitude equation on its natural time scale
in case of $ d $-dimensional  periodic domains of length 
$ \mathcal{O}(1/\varepsilon) $ for the Swift-Hohenberg model under suitable conditions 
on the additive  noise.
In detail, we prove the validity of this approximation for noise 
whose set of Fourier coefficients with respect to $ x $ is in $ \ell^1 $ 
for fixed $ t \geq 0 $. Moreover, we improve existing approximation results
in the sense that the stable part of the noise can be  larger.
\end{abstract}


\section{Introduction}

We start with an anisotropic $ d $-dimensional Swift-Hohenberg model 
\begin{equation} \label{SH}
\partial_t u = - (1+\partial_{x_1}^2)^2 u + 4 \Delta_{x_{\perp}} u + \alpha u -u^3 + \zeta 	,
\end{equation}
with $ x = (x_1, x_{\perp} ) = (x_1,\ldots, x_d) \in \R^d $,  $ u(x,t) \in \R $, 
$ \Delta_{x_{\perp}} = \sum_{j=2}^d \partial_{x_j}^2 $, $ t \geq 0 $, and additive  noise $ \zeta = \zeta(x,t) $
in time and space,
which will be  defined by a Fourier series in space and whose coefficients are distributional derivatives of Wiener processes in time. 
The noise will be  specified more precisely below. Moreover, $ \alpha \in \R $ is a control parameter.

For vanishing noise $ \zeta = 0 $
the linearized problem  at  $ u = 0 $ is solved 
by $$ u(x,t) = e^{\lambda(k) t +
\alpha t + ik \cdot x} ,$$ with $ k = (k_1,k_{\perp}) = (k_1,\ldots,k_d) $,
$ k \cdot x = \sum_{j = 1}^d k_j x_j $, and
$$ \lambda(k) = -(1-k_1^2)^2 - 4 |k_{\perp}|^2 .$$
Therefore, the trivial solution $ u = 0 $ becomes unstable if $ \lambda(k) + \alpha  > 0 $ for a 
$ k \in \R $. We find that $  \lambda(k) = 0 $ 
at the wave vectors $ (k_1,k_{\perp}) = (\pm 1,0)$. 
This kind of instability appears for instance in electro-convection for which $ d = 2 $, cf. \cite{OD08}.
In order 
to describe the bifurcating solutions 
for small $ \alpha > 0 $ we 
introduce the small perturbation parameter $ 0 < \varepsilon \ll 1 $, set $ \alpha = \varepsilon^2 $ 
and 
make the ansatz
\begin{equation}\label{ansatz}
u(x,t) = \varepsilon A(X,T) e^{ix_1}
+ c.c. ,
\end{equation}
with slow spatial variables $ X = \varepsilon x$,
slow time variable $ T = \varepsilon^2 t $, and complex valued amplitude 
$ A(X,T) \in \C $.
We find that in lowest order $ A $
has to satisfy a stochastic Ginzburg-Landau equation 
\begin{equation}\label{GL}
	\partial_T A = 4 \Delta_X A  + A - 3 A |A|^2 +
	\zeta_A,
\end{equation}
with $ \Delta_X = \sum_{j=1}^d \partial_{X_j}^2 $ and $ \zeta_A= \zeta_A(X,T)  $ being a rescaled 
noise process in time and space which will be also specified below.

We consider 
\eqref{SH} with $  2 \pi/ \varepsilon $-periodic boundary conditions
\begin{equation} \label{bcSH}
u(x,t) = u(x+ 2 \pi e_j/ \varepsilon,t) 
\end{equation}
for $ j = 1,\ldots,d $,
where $ e_j $ is the $ j $-th unit vector.
These boundary conditions correspond to $ 2 \pi $-periodic  boundary conditions
\begin{equation} \label{bcGL}
A(X,T) = A(X+ 2 \pi e_j,T)  
\end{equation}
for $ j = 1,\ldots,d $
for the stochastic Ginzburg-Landau equation \eqref{GL}. 
In contrast to the whole of $ \R^d $ these large 
periodic domains allow us  to define  and to control the noise 
as Fourier series in space whose coefficients are white noise processes in time. Our approach is fundamentally based on 
this representation of the noise as Fourier series exactly like in the derivation and justification of the 
Ginzburg-Landau approximation which is also fundamentally based on 
the representation of the solutions in Fourier space.
These kind of boundary conditions have been used before, cf. \cite{MSZ00} and also in \cite{BHP05} where 
the 1D case has been handled.
It will be the topic of future research how the analysis from \cite{BB19,BBS19} which is not based 
on Fourier transform can be transferred to the handling of \eqref{SH}
on $ \R^d $.

It is the goal of this paper to prove that the stochastic Ginzburg-Landau equation \eqref{GL} with the boundary conditions 
\eqref{bcGL} makes correct predictions about the dynamics of the stochastic Swift-Hohenberg model \eqref{SH}
with the boundary conditions 
\eqref{bcSH}.  Already in one space-dimension the validity of the stochastic 
Ginzburg-Landau approximation turns out to be a subtle problem
due to the low regularity of the involved stochastic processes in time and space, cf.
\cite{BHP05}.
This problem becomes worse in higher space dimensions.
In the following we would like to present an approach which can be generalized 
to arbitrary space dimensions. Assuming 
that the noise 
is in some Wiener algebra for fixed $ t \geq 0 $, i.e., assuming 
that the 
set of  Fourier coefficients is in $ \ell^1 $ 
for fixed $ t \geq 0 $ allows us to transfer the existing approximation theory from 
the deterministic situation to the stochastic situation.

By working with projections which separate the exponentially damped from the weakly 
unstable modes we can show that the noise in the stable part can 
be of order $ \mathcal{O}(\varepsilon^{\beta}) $ for any fixed $ \beta > 1 $.
This improves the result from \cite{BHP05} where $ \mathcal{O}(\varepsilon^{5/2}) $ in $ L^2 $ was assumed for the noise.
This can be interpreted as an uncertainty quantification result.

The plan of the paper is as follows.
In Section \ref{sec2} we recall the derivation of the Ginzburg-Landau equation in the deterministic case.
In Section \ref{sec3} we discuss
the properties of the solutions of the stochastic Ginzburg-Landau equation. We introduce
the used function spaces, we estimate
the appearing stochastic terms, and present the
local existence and uniqueness theory for solutions of this stochastic PDE.
In Section \ref{sec4} we do some preparations for the Swift-Hohenberg equation, like 
the ones we do in Section \ref{sec3} for the stochastic Ginzburg-Landau equation \eqref{GL}.
Section \ref{sec5} contains the approximation result and its proof.
The paper is closed with an outlook and discussion section.
\medskip

{\bf Notation:} Throughout this paper many possibly different constants 
will be denoted with the same symbol $C$ if they can be chosen independently 
of the small perturbation parameter $ 0 < \varepsilon \ll 1 $.
In the following we work with functions which satisfy either the boundary conditions \eqref{bcGL}
or \eqref{bcSH}. In most cases
functions satisfying \eqref{bcGL} are denoted with large capitals such as $ A, B, \ldots $.
In most cases functions satisfying \eqref{bcSH} are denoted with small capitals such as $ u, v, \ldots $.
Their Fourier transforms w.r.t. \eqref{bcGL} are denoted with $ \widehat{A}, \widehat{B}, \ldots $.
Their Fourier transforms w.r.t. \eqref{bcSH} are denoted with $ \widehat{u}, \widehat{v}, \ldots $.
\medskip

{\bf Acknowledgement.}
The authors are grateful for discussions with Dirk Bl\"omker and Reika Fukuizumi.
The work of 
Guido Schneider is partially supported by the Deutsche Forschungsgemeinschaft DFG through the cluster of excellence ’SimTech’ under EXC 2075-390740016.

\section{Derivation of the Ginzburg-Landau equation in the deterministic case}
\label{sec2}

Before we answer the question posed in the introduction we recall the derivation of the Ginzburg-Landau  equation in the deterministic case.
For the deterministic anisotropic Swift-Hohenberg equation 
\begin{equation} \label{SHdet}
\partial_t u = - (1+ \partial_{x_1}^2)^2 u + 4 \Delta_{x_{\perp}} u + \varepsilon^2 u - u^3
\end{equation}
we define the residual
$$
\textrm{Res}_{det}(u) = - \partial_t u - (1+ \partial_{x_1}^2)^2 u + 4 \Delta_{x_{\perp}} u + \varepsilon^2 u - u^3
$$
which counts how much a function $ u $ fails to satisfy the Swift-Hohenberg equation \eqref{SHdet}. 
For the derivation of the  Ginzburg-Landau equation
we make the ansatz 
\begin{eqnarray*}
u_A(x,t) & = & \varepsilon A_1(X,T) e^{ix_1} + \varepsilon A_{-1}(X,T) e^{-ix_1},
\end{eqnarray*}	
where  $ X = \varepsilon x $, $ T = \varepsilon^2 t $,
and where $ A_{-1} = \overline{A}_1 $. 
\begin{remark}{\rm
If $ A_1 $ satisfies the boundary conditions \eqref{bcGL}, then $ u_A $ satisfies the boundary conditions \eqref{bcSH}
only if $ \varepsilon = 1/n $ for $ n \in \N $. To avoid this restriction we can replace \eqref{bcSH} by 
\begin{equation} \label{bcSHp}
u(x,t) = u(x+ 2 \pi e_j/ \varepsilon + b(\varepsilon),t) 
\end{equation}
with $ b(\varepsilon) \in \R^d $ suitably chosen satsifying $ \limsup_{\varepsilon \to 0} |b(\varepsilon)| = \mathcal{O}(1) $.
For notational simplicity the subsequent analysis is carried out for \eqref{bcSH}. The generalization to \eqref{bcSHp}
is obvious.
}
\end{remark}
We find
\begin{eqnarray*}
\textrm{Res}_{det}(u_A) & = & - \partial_t u_A -  u_A - 2 \partial_{x_1}^2 u_A - \partial_{x_1}^4 u_A 
 + 4 \Delta_{x_{\perp}} u_A + \varepsilon^2 u_A - u_A^3 \\
& = &  s_1+ \ldots + s_{7} ,
\end{eqnarray*}	
where
\begin{eqnarray*}
s_1 & = &  - \varepsilon^3 \partial_T  A_1 e^{ix_1} -  \varepsilon^3 \partial_T A_{-1} e^{-ix_1}, 
\\
s_2 & = & 
- \varepsilon A_1 e^{ix_1}-  \varepsilon A_{-1} e^{-ix_1} ,
\\ 
s_3 & = &
  2 \varepsilon   A_1 e^{ix_1} + 2 \varepsilon   A_{-1} e^{-ix_1} 
 - 4 i  \varepsilon^2  \partial_{X_1}  A_1 e^{ix_1} + 4 i \varepsilon^2 \partial_{X_1} A_{-1} e^{-ix_1} 
\\&& 
 - 2 \varepsilon^3 \partial_{X_1}^2  A_1 e^{ix_1} - 2 \varepsilon^3 \partial_{X_1}^2 A_{-1} e^{-ix_1} ,
\\
s_4 & = &  -  \varepsilon   A_1 e^{ix_1} -  \varepsilon   A_{-1} e^{-ix_1} 
+ 4 i  \varepsilon^2  \partial_{X_1}  A_1 e^{ix_1} - 4 i \varepsilon^2 \partial_{X_1} A_{-1} e^{-ix_1} 
\\&& 
+ 6 \varepsilon^3 \partial_{X_1}^2  A_1 e^{ix_1} + 6 \varepsilon^3 \partial_{X_1}^2 A_{-1} e^{-ix_1} 
+ 4 i \varepsilon^4 \partial_{X_1}^3  A_1 e^{ix_1} - 4 i \varepsilon^4 \partial_{X_1}^3 A_{-1} e^{-ix} 
\\&& 
-  \varepsilon^5 \partial_{X_1}^4  A_1 e^{ix_1} -  \varepsilon^5 \partial_{X_1}^4 A_{-1} e^{-ix_1} ,
\\ 
s_5 & = & 
4 \varepsilon^3 \Delta_{X_{\perp}}  A_1 e^{ix_1} + 4\varepsilon^3 \Delta_{X_{\perp}} A_{-1} e^{-ix_1}, 
\\
s_6 & = & -  \varepsilon^2 s_2,
\\
s_{7} & = & -( \varepsilon A_1 e^{ix_1}+  \varepsilon A_{-1} e^{-ix_1} )^3.
\end{eqnarray*}	
We see that the terms of order $ \mathcal{O}(\varepsilon) $ from $ s_2 $, $ s_3 $, and $ s_4 $
cancel.
The same happens to the terms of order $ \mathcal{O}(\varepsilon^2) $ from
$ s_3 $ and $ s_4 $.
The terms of order $ \mathcal{O}(\varepsilon^3) $ from $ s_1,s_3,\ldots,s_7 $ with a factor $ e^{ix_1} $ give
 the  Ginzburg-Landau equation 
\begin{equation} \label{GLdet}
0 = - \partial_T  A_1 + 4  \Delta_X  A_1   + A_1 - 3 A_1 |A_1|^2 
\end{equation}
by equating their sum to zero. 
We obtain the complex conjugate equation for $ A_{-1} $ by equating the terms of order $ \mathcal{O}(\varepsilon^3) $  with a factor $ e^{-ix} $ to zero. 
Hence,  we finally have 
\begin{eqnarray*}
\textrm{Res}_{det}(u_A) & = &  4 i \varepsilon^4 \partial_{X_1}^3  A_1 e^{ix_1} - 4 i \varepsilon^4 \partial_{X_1}^3 A_{-1} e^{-ix_1} \\ &&
-  \varepsilon^5 \partial_{X_1}^4  A_1 e^{ix_1} -  \varepsilon^5 \partial_{X_1}^4 A_{-1} e^{-ix_1} 
\\ && - \varepsilon^3 A_1^3 e^{3ix_1}-  \varepsilon^3 A_{-1}^3 e^{-3ix_1}  .
\end{eqnarray*}	
\begin{remark}
\label{rem21}
{\rm 
In the deterministic case 
the validity of the Ginzburg-Landau approximation with an error of order $ \mathcal{O}(\varepsilon^2) $ 
follows with a simple application of Gronwall's inequality by considering an extended ansatz 
\begin{eqnarray*}
u_A(x,t) & = & \varepsilon A_1(X,T) e^{ix_1} + \varepsilon A_{-1}(X,T) e^{-ix_1},
\\ &&
+  \varepsilon^3 A_3(X,T) e^{3ix_1} + \varepsilon^3 A_{-3}(X,T) e^{-3ix_1},
\end{eqnarray*}	
with $ 0 = - 64 A_3 + A_1^3 $. Then in  the residual only terms  of order $ \mathcal{O}(\varepsilon^4) $ remain.
However,   the term $ \partial_T A_3 $ has  to 
be estimated. If $ A_1 $ is sufficiently regular in space, this is not a problem. However, in the stochastic case 
the necessary regularity is not present and so we refrain from using the extended ansatz but use an approach with mode
filters instead.
}\end{remark}


%
%
%

\section{The stochastic Ginzburg-Landau equation}
\label{sec3}

We extend the   Ginzburg-Landau equation 
with a 
spatio-temporal noise term $ \zeta_A = \zeta_A(X,T) $ which is specified in the following.
Therefore, we consider 
\begin{equation} \label{2DGL}
 \partial_T  A_1 =  4  \Delta_X  A_1  + A_1 - 3 A_1 |A_1|^2 + \zeta_A,
\end{equation}
where $ X  \in \R^d $, $ T \geq 0 $, and $ A_1(X,T) \in \C $.
This equation is equipped with $ 2 \pi $-periodic  boundary conditions 
$$
A_1(X,T) = A_1(X+ 2 \pi e_j,T)  ,
$$
for $ j = 1, \ldots,d $.
Similarly, we assume 
$$
\zeta_A(X,T) = \zeta_A(X+ 2 \pi e_j,T) 
$$
for $ j = 1, \ldots,d $ for the noise term.
Due to the  
periodicity of the involved functions we use Fourier series
\begin{eqnarray*}
 A_1(X,T)  & = & \sum_{K \in \Z^d} 
 \widehat{A}_1(K,T) e^{iK \cdot X} , \\
 \zeta_A(X,T)  & = & \sum_{K \in \Z^d} 
 \widehat{\zeta}_A(K,T) e^{iK \cdot X} .
\end{eqnarray*}

\subsection{The function spaces}

The magnitude of these functions will be controlled 
in Fourier space 
by $ \ell^1_{r} $-norms defined by
$$
\|  \widehat{B}\|_{\ell^1_r}
:= \sum_{K \in \Z^d} | \widehat{B}(K) | (1+|K|^2)^{r/2}.
$$
for $ r \in \R $.
We use the abbreviation $ \ell^1 = \ell^1_0 $ and
sometimes we use the notation
$$
\|\widehat{B}\|_{L^1_{A,r}} = \int_{\mathbb{R}^d}
| \widehat{B}(K) | (1+|K|^2)^{r/2}d\mu_A(K)  
:= \sum_{K \in \Z^d} | \widehat{B}(K) | (1+|K|^2)^{r/2}
$$
to make more transparent  in the following the similarity of our approach to the pure deterministic case. Herein, $ d\mu_A(k) = \sum_{k' \in \Z^d} \delta(k-k') $ is a measure with point masses $ 1 $
on the lattice $ \Z^d $.


Fourier transform is continuous from $ \ell^1_r $ to $ C^r_b $ for $ r \in \N $,
but not vice versa.
The space $ \ell^1_r $ is closed 
under discrete convolution  
$$
(\widehat{u} *\widehat{v} )(k) = \sum_{l \in \Z^d} \widehat{u}(k-l) \widehat{v} (l) 
$$
if $ r \geq 0$. This is based on 
Young's inequality
$$
\| \widehat{u} *\widehat{v}  \|_{\ell^1} \leq C  \| \widehat{u}   \|_{\ell^1}\|  \widehat{v}  \|_{\ell^1},
$$
and on the inequality
$$ 
(1 +k^2)^r \leq C ((1 +(k-l)^2)^r  + (1 +l^2)^r ) 
$$ 
for $ r \geq 0 $. In detail, we have
that 
for $  r \geq 0 $ 
 there 
exists a $ C > 0 $ such that for all $ \widehat{u},\widehat{v} \in \ell^1_r $ 
we have 
$$ 
\| \widehat{u} *\widehat{v}  \|_{\ell^1_r} \leq C \| \widehat{u}   \|_{\ell^1_r}\|  \widehat{v}  \|_{\ell^1_r}.
$$ 
The counterpart in physical space are the Banach spaces 
$ \mathcal{W}_A^r = \mathcal{W}_A^r(\mathbb{T}^d)  $ equipped with the norm 
$$ 
\|B \|_{\mathcal{W}_A^r} = \|\widehat{B}\|_{L^1_{A,r}} .
$$ 
In the literature,  these spaces sometimes are called Wiener algebras.
For $  r \geq 0 $ it is closed under point-wise multiplication as a direct consequence of the above estimates, i.e., 
for all $ r \geq 0 $ there 
exists a $ C > 0 $ such that for all $ u,v \in \mathcal{W}_A^r $ 
we have 
$$ 
\| u v \|_{\mathcal{W}_A^r} \leq C \| u  \|_{\mathcal{W}_A^r}\| v  \|_{\mathcal{W}_A^r}.
$$

\subsection{Estimates for the stochastic terms}

To rigorously define the  $ 2 \pi $-spatially periodic noise term $ \zeta_A = \zeta_A(X,T) $
we split it in its real and imaginary part, namely 
\begin{eqnarray*}
\widehat{\zeta}_A(K,T)  &=& \alpha_A(K)  \widehat{\xi}_A(K,T) = 
\widehat{\zeta}_{A,r}(K,T) +  i
\widehat{\zeta}_{A,i}(K,T) \\ & = & \alpha_A(K) (\widehat{\xi}_{A,r}(K,T) +  i
\widehat{\xi}_{A,i}(K,T) ),
\end{eqnarray*}
where  $ \widehat{\xi}_{A,r}(K,\cdot) $ and $ \widehat{\xi}_{A,i}(K,\cdot) $
are distributional derivatives w.r.t. to time $ T \geq 0 $ of standard Wiener processes 
$ \widehat{W}_{A,r}(K,\cdot) $ and $ \widehat{W}_{A,i}(K,\cdot) $
for each $ K \in \Z^d  $. The coefficients 
$  \alpha_A(K) \geq  0 $  
are specified below and could allow us for instance to control firstly the spatial regularity and secondly the magnitude with respect to 
$ \varepsilon $ of the noise process.

For solving the stochastic Ginzburg-Landau equation we separate the lowest order stochastic part from the solution by
introducing $ A_1 = B + Z_A $, with $ B $ and $ Z_A $ satisfying 
\begin{eqnarray} \label{eqBneu}
\partial_T B & = & 4 \Delta_X B + B + 2Z_A- 3 |B+Z_A|^2 (B+Z_A)   , \\ 
\partial_T Z_A & = & 4 \Delta_X Z_A - Z_A + \zeta_A  , \label{eqZAneu}
\end{eqnarray}
The first equation will be handled in the next section.
Estimates for the second equation will be provided here.

In Fourier space the second equation is given by
$$ 
\partial_T \hZ_A(K,T)    =  - 4|K|^2 \hZ_A(K,T) -  \hZ_A(K,T) +  \alpha_A(K)\widehat{\xi}_A(K,T)
$$ 
which  will be solved with the help of  the variation of constant formula.
Using $  \widehat{Z}_A(K,0) = 0 $ we find for $ K \in \Z^d  $ that
\begin{eqnarray} \nonumber
\widehat{Z}_A(K,T) & = & \int_0^T  e^{-(1+ 4 |K|^2) (T-\tau)}  \alpha_A(K) \widehat{\xi}_A(K,\tau) d\tau
\\ & = & 
\int_0^T  e^{-(1+ 4 |K|^2) (T-\tau)} \alpha_A(K) d\widehat{W}_A(K,\tau)  , \label{Zeqvo}
\end{eqnarray}
where the integral is a stochastic integral in the It\^{o}-sense.
For fixed $ K \in \Z^d $ the process
$ (\widehat{Z}_A(K,T))_{T \geq 0} $ is a so called  Ornstein-Uhlenbeck
process. These are well studied stochastic processes for which 
we recall some important estimates from the existing literature.

We assume 
\begin{equation} \label{ASS1}
 (\alpha_A(K))_{K \in \Z^d} \in \ell^1_{r_A} 
\end{equation}
where $ r_A = 0 $ will be chosen below.
In the following 
we need estimates for 
$$ 
\sup_{\tau \in [0,T]}\|\widehat{Z}_A(\cdot,\tau) \|_{\ell^1_r} ,
$$ 
with $ r \geq 0 $ as large as possible. 
Doob's submartingal theorem, cf. \cite[Theorem 1.51]{Baudoin}, combined with \cite[Exercise 1.33]{Baudoin} yields
$$ 
\text{P}[\sup_{\tau \in [0,T]}\|\widehat{Z}_A(\cdot,\tau) \|_{\ell^1_r} \geq c]\leq \frac{\mathbb{E}[\|\widehat{Z}_A(\cdot,T) \|_{\ell^1_r}^2]}{c^2}= s_1,
$$
with
$$
\mathbb{E}[u] = \int_{\Omega} u(\omega) dP_A(\omega) 
$$
where $ \Omega $ is the space of realizations of the noise process 
with associated probability measure $ P_A(\omega) $, cf. \cite{Prato}.

Using H\"older's inequality in the third estimate and the assumption \eqref{ASS1} in the fourth estimate,
for $ s_1 $ we estimate 
 \begin{eqnarray*}
s_1 & \leq & \frac{1}{c^2} \int |\int_{\R^d} |\widehat{Z}_A(K,T,\omega)| (1+|K|^2)^{r/2}  d\mu_A(K) |^2 dP_A(\omega)  \\
& = & \frac{1}{c^2} \int |\int_{\R^d} |\widehat{Z}_A(K,T,\omega)| (1+|K|^2)^{r/2} \alpha_A(K)^{-1/2} \\ && \qquad \qquad  \qquad  \times \alpha_A(K)^{1/2} d\mu_A(K) |^2 dP_A(\omega)  \\
& \leq & \frac{1}{c^2} \int |\int_{\R^d}  |\widehat{Z}_A(K,T,\omega)|^2 (1+|K|^2)^{r-r_A/2}  \alpha_A(K)^{-1}d\mu_A(K)  \\ && \qquad \qquad \times \int_{\R^d}  \alpha_A(K) (1+|K|^2)^{r_A/2}  d\mu_A(K) | dP_A(\omega)  \\
& \leq  & C \frac{1}{c^2} \int |\int_{\R^d}  |\widehat{Z}_A(K,T,\omega)|^2 (1+|K|^2)^{r-r_A/2} \alpha_A(K)^{-1} d\mu_A(K)   | dP_A(\omega)  \\
& \leq  & C \frac{1}{c^2} \int_{\R^d}  \int |\widehat{Z}_A(K,T,\omega)|^2  dP_A(\omega)  (1+|K|^2)^{r-r_A/2}    \alpha_A(K)^{-1}   d\mu_A(K) \\
& =  & C \frac{1}{c^2} \int_{\R^d}  \mathbb{E}|\widehat{Z}_A(K,T,\omega)|^2 (1+|K|^2)^{r-r_A/2}    \alpha_A(K)^{-1}  d\mu_A(K) = s_2.
\end{eqnarray*}
Applying It\^{o}'s isometry 
$$
\mathbb{E} \left[ \left( \int_0^T X_t \, \mathrm{d} W_t \right)^2 \right] 
= \mathbb{E}\left[ \int_0^T X_t^2 \, \mathrm{d} t \right]
$$
to $ s_2 $, using \eqref{Zeqvo}, gives 
 \begin{eqnarray*}
c^2 s_2 & = & C \int_{\R^d}  \mathbb{E} \left[ \left( \int_0^T e^{-(1+ 4 |K|^2) (T-\tau)} \alpha_A(K)  d\widehat{W}(K,\tau)   \right)^2 \right] \\ && \qquad \qquad \times(1+|K|^2)^{r-r_A/2} \alpha_A(K)^{-1}  d\mu_A(K) \\
& = & C \int_{\R^d}  \mathbb{E} \left[ \left( \int_0^T e^{-(1+ 4 |K|^2)(T-\tau)} d\widehat{W}(K,\tau)   \right)^2 \right] \\ && \qquad \qquad \times(1+|K|^2)^{r-r_A/2} \alpha_A(K)  d\mu_A(K) \\
& = & C \int_{\R^d}   \left[ \int_0^T (e^{-(1+ 4 |K|^2) (T-\tau)})^2 \, d \tau \right] (1+|K|^2)^{r-r_A/2} \alpha_A(K)  d\mu_A(K)
\\ 
& = & C \int_{\R^d} \left[ \int_0^T (e^{-2 (1+ 4 |K|^2) (T-\tau)}  \, d \tau \right]  \alpha_A(K)  (1+|K|^2)^{r-r_A/2}  d\mu_A(K)
\\ & \leq & 
C \int_{\R^d} \frac{\alpha_A(K) (1+|K|^2)^{r - r_A/2}}{(1+ 4 |K|^2)} d\mu_A(K)
\\ & \leq & 
C \sup_{K \in \R^d} \frac{(1+|K|^2)^{r-r_A}}{(1+ 4 |K|^2)} \times
\int_{\R^d} \alpha_A(K) (1+|K|^2)^{r_A/2}  d\mu_A(K)
 < \infty
\end{eqnarray*}
if $$ r-r_A \leq 1, $$
respectively 
$$ r \leq  r_A +1 ,$$
again due to \eqref{ASS1}.
We  finally get the lemma.
\begin{lemma}
Suppose $ r \leq  r_A +1 $.
Then for all $ \delta_0 > 0 $ there 
exists a $ c > 0 $ such that 
$$
\text{P}[\sup_{\tau \in [0,T]}\|\widehat{Z}_A(\cdot,\tau) \|_{L^1_{A,r}} \geq c] \leq \delta_0,
$$
respectively 
$$
\text{P}[\sup_{\tau \in [0,T]}\|Z_A(\cdot,\tau) \|_{\mathcal{W}_{A}^r} \geq c] \leq \delta_0.
$$
\end{lemma}
\begin{corollary} \label{cor34a}
For  $ r_A = 0 $ and all $ \delta_0 > 0 $ there 
exists a $ C_{Z_A} > 0 $ and a $ T_0 > 0 $ such that 
we have solutions $ Z_A  $ 
of \eqref{Zeqvo} with 
$$
\text{P}[\sup_{\tau \in [0,T_0]}\|Z_A(\cdot,\tau) \|_{\mathcal{W}_A^{1}} \leq C_{Z_A}] \geq 1- \delta_0.
$$
\end{corollary}
Note that we have the continuity in time of $ Z_A(\cdot,t) $ in the 
space $ \mathcal{W}_A^{1} $. 
This is obviously true for any finite-dimensional approximation of $ Z_A $ 
and therefore also for $ Z_A $ as a uniform limit in time of these finite-dimensional approximations.
Hence, the subsequent integrations 
in the variation of constant formulas w.r.t. time are well-defined.


\subsection{Local existence and uniqueness}

With $ A_1 = B + Z_A $ we have rewritten  the stochastic Ginzburg-Landau equation 
$$
 \partial_T  A_1 =  4  \Delta  A_1  + A_1 - 3 A_1 |A_1|^2 + \zeta_A 
$$ 
as
\begin{eqnarray} 
\partial_T B & = & 4 \Delta B + B + 2Z_A- 3 |B+Z_A|^2 (B+Z_A)   , \\ 
\partial_T Z_A   & = & 4 \Delta Z_A - Z_A + \zeta_A  . 
\end{eqnarray}
For handling the first equation we consider the variation of constant formula associated to \eqref{eqBneu}, in Fourier space, namely 
\begin{eqnarray} \label{vocnewGL}
\lefteqn{\widehat{B}(K,T)  =  e^{- (4 | K |^2 + 1) T} \widehat{B}(K,0)} \\ && + 
\int_0^T e^{- (4 | K |^2 + 1)(T-\tau)}
(2\widehat{Z}_A -3  (\widehat{B}+ \widehat{Z}_A) * (\widehat{B} + \widehat{Z}_A)
* (\widehat{\overline{B}}+ \widehat{\overline{Z}}_A))(K,\tau) d\tau.  \nonumber
\end{eqnarray}

For the construction of the solution we first use that 
the space $ \ell^1_r $ is closed under convolution for $ r \geq 0 $. 
Secondly, using the smoothing property of the linear semigroup $ (e^{- (4 | K |^2 + 1) T})_{T \geq 0} $ from $ \ell^1_r $ to $ \ell^1_{r+2-\delta'} $ for a fixed but arbitrary small $ \delta' > 0 $,
with a singularity 
$ (t-\tau)^{-(2-\delta')/2} $, shows that the right-hand side of \eqref{vocnewGL}
is a contraction in a ball in 
$C([0,T_0],\ell^1_{r+2-\delta'}) $
if $ T_0 = T_0(C_{Z_A},\| \widehat{B}(\cdot,0) \|_{\ell^1_{r+2-\delta'}} ) > 0 $ is sufficiently small.

Therefore, we have 
\begin{theorem}
Fix $ r = r_A + 1 > 0 $. 
Then for a fixed but arbitrary small $ \delta' > 0 $ and  for all $ \delta_0 > 0 $ there 
exists a $ c > 0 $ and a $ T_0 > 0 $ such that 
we have solutions $ \widehat{B} $ 
of \eqref{vocnewGL} with 
$$
\text{P}[\sup_{\tau \in [0,T_0]}\|\widehat{B}(\cdot,\tau) \|_{\ell^1_{r+2-\delta'}} \leq c] \geq 1- \delta_0.
$$
\end{theorem}
\begin{corollary} \label{cor34}
For  $ r_A = 0 $, for a fixed but arbitrary small $ \delta' > 0 $, and all $ \delta_0 > 0 $ there 
exists a $ C_B > 0 $ and a $ T_0 > 0 $ such that 
we have solutions $ B $ 
of \eqref{vocnewGL} with 
$$
\text{P}[\sup_{\tau \in [0,T_0]}\|{B}(\cdot,\tau) \|_{\mathcal{W}_A^{3-\delta'}} \leq C_B] \geq 1- \delta_0.
$$
\end{corollary}

\section{The Wiener algebra}

\label{sec4}

It is essential to have a suitable functional analytic set-up for handling 
the periodic boundary conditions \eqref{bcSH}.  It turns out that the situation 
of a domain of size $ \mathcal{O}(1/\varepsilon) $
is more close to the unbounded domain than to the situation of a domain of size 
$ \mathcal{O}(1) $.  

\subsection{Definition and basic properties}

For this reason and 
in order to have the correct scalings
w.r.t. $ \varepsilon $ we base our analysis on the 
continuous inverse Fourier transform  
$$ 
u (x)= \int_{\R^d} \widehat{u}(k) e^{i k\cdot x}
d\mu(k)
$$
and the Fourier transform  
$$ 
\widehat{u}(k)   = \frac{1}{(2\pi)^d} \int_{\mathbb{T}^d_{2 \pi/\varepsilon}} u(x)  e^{-i k\cdot x} d x,
$$
where with $\mathbb{T}_{2 \pi/\varepsilon} = \R/(2 \pi \Z/\varepsilon)  $
we took care of the periodic boundary conditions \eqref{bcSH}.
Moreover, $ d\mu(k) $ is a measure with point masses $ \varepsilon^d $
on the lattice 
$$
\mathcal{L}_{\varepsilon}
= \{  k = \varepsilon j_k:  j_k \in \Z^d \},
$$
see below for more details.
Due to length of the domain the 
Fourier coefficients of a function which is $ \mathcal{O}(1) $ in physical space can have
Fourier coefficients of order $ \mathcal{O}(\varepsilon^{-d})$.
Therefore, it is advantageous to stay as close as possible to the analysis 
used on unbounded domains and so 
we write 
the inverse Fourier transform 
as 
\begin{eqnarray*}
 u(x)  & = & \int_{\R^d} \widehat{u}(k) e^{ik \cdot x} d\mu(k) := 
 \sum_{j_k \in \Z^d} 
 \widehat{u}(\varepsilon j_k) e^{i \varepsilon j_k \cdot x}\Delta k 
 \\ & := & \varepsilon^d \sum_{j_k \in \Z^d} 
 \widehat{u}(\varepsilon j_k)e^{i \varepsilon j_k \cdot x}
 = \varepsilon^d \sum_{k \in \mathcal{L}_{\varepsilon}} 
 \widehat{u}(k) e^{ik \cdot x} .
\end{eqnarray*}
For the same reasons
our function spaces should be more close to $ C^0_b $ 
than to $ L^2 $.  
The counterpart to $ C^0_b$ in 
physical space is $ L^1 $ in Fourier space, here equipped with the norm
$$
\|  \widehat{u}\|_{L^1} = \int_{\R^d}
| \widehat{u}(k) | d\mu(k) 
:= \varepsilon^d\sum_{k \in \mathcal{L}_{\varepsilon}} | \widehat{u}(k) |.
$$
We have the estimate 
\begin{eqnarray*}
\sup_{x \in \R^d}|u(x)|  & = & \sup_{x \in \R^d} | \int_{\R^d} \widehat{u}(k) e^{ik \cdot x} d\mu(k) | \\
& \leq & \varepsilon^d \sum_{k \in \mathcal{L}_{\varepsilon}} 
 \sup_{x \in \R^d} |  \widehat{u}(k) e^{ik\cdot x} |
\leq  \varepsilon^d \sum_{k \in \mathcal{L}_{\varepsilon}}  |  \widehat{u}(k) | \\
& = &  \int_{\R^d} |\widehat{u}(k)| d\mu(k).
\end{eqnarray*}
Hence, inverse Fourier transform works like on the real line and is continuous from $ L^1 $ to $ C^0_b $ with an $ \mathcal{O}(1) $-bound
but not vice versa. Like for the Ginzburg-Landau equation we define the space $ L^1_r $ equipped with the norm 
$$ 
\|  \widehat{u}\|_{L^1_r} = \int_{\R^d}
| \widehat{u}(k) | (1+|k|^2)^{r/2} d\mu(k) 
= \varepsilon^d\sum_{k \in \mathcal{L}_{\varepsilon}} | \widehat{u}(k) | (1+|k|^2)^{r/2}.
$$ 
Like on the real line the inverse Fourier transform  is continuous from $ L^1_r $ to $ C^r_b $ with an $ \mathcal{O}(1) $-bound for $ r \in \N_0 $.
We define the Banach space $ \mathcal{W}^r   $ equipped with the norm 
$$ 
\| u \|_{\mathcal{W}^r} = \|\widehat{u}\|_{L^1_{r}} .
$$ 
Like for the Ginzburg-Landau equation we have that 
for all $ r \geq 0 $ there 
exists a $ C > 0 $ such that for all $ u,v \in \mathcal{W}^r $ 
we have 
$$ 
\| u v \|_{\mathcal{W}^r} \leq C \| u  \|_{\mathcal{W}^r}\| v  \|_{\mathcal{W}^r}.
$$

\subsection{Scaling properties}

The Fourier transform of the ansatz 
\begin{eqnarray*}
u_A(x,t) & = & \varepsilon A_1(\varepsilon x,\varepsilon^2 t) e^{ix_1} + \varepsilon A_{-1}(\varepsilon x,\varepsilon^2 t) e^{-ix_1}
\end{eqnarray*}	
is given by 
$$
\widehat{u}_A(k,t) =  \varepsilon \varepsilon^{-d} \widehat{A}_1(\frac{k-e_1}{\varepsilon},\varepsilon^2 t) + \varepsilon  \varepsilon^{-d} \widehat{A}_{-1}(\frac{k+e_1}{\varepsilon},\varepsilon^2 t) .
$$
As a direct consequence of this scaling behavior and the definition of our norms we have for $ r = 0 $ that 
\begin{eqnarray*}
\| u_A(\cdot,t) \|_{\mathcal{W}} & = &  \int_{\R^d}
| \widehat{u}_A(k,t) | d\mu(k) 
= \varepsilon^d\sum_{k \in \mathcal{L}_{\varepsilon}} | \widehat{u}_A(k,t) |
\\ 
& = & \varepsilon^d\sum_{k \in \mathcal{L}_{\varepsilon}} |
\varepsilon \varepsilon^{-d} \widehat{A}_1(\frac{k-e_1}{\varepsilon},\varepsilon^2 t)+ \varepsilon  \varepsilon^{-d} \widehat{A}_{-1}(\frac{k+e_1}{\varepsilon},\varepsilon^2 t) | \\
& \leq & 2 \varepsilon^d\sum_{k \in \mathcal{L}_{\varepsilon}} |
\varepsilon \varepsilon^{-d} \widehat{A}_1(\frac{k-e_1}{\varepsilon},\varepsilon^2 t) |
= 2 \sum_{k \in \mathcal{L}_{\varepsilon}} |
\varepsilon  \widehat{A}_1(\frac{k-e_1}{\varepsilon},\varepsilon^2 t) |
\\ & = & 2 \sum_{K \in \Z^d} |
\varepsilon  \widehat{A}_1(K,\varepsilon^2 t) |
= 2 \varepsilon  \|{A}_1 \|_{\mathcal{W}_A},
\end{eqnarray*}	
respectively
$$ 
\| u_A(\cdot,t) \|_{\mathcal{W}^r} \leq C  \varepsilon  \| {A}_1 \|_{\mathcal{W}_A^r}
$$ 
for all $ r \geq 0 $ with a constant only depending on $ r $.

\section{The approximation result}
\label{sec5}

This section contains the proof of our approximation result.

\subsection{Connection between $ \zeta $ and $ \zeta_A $}

\paragraph{The projections.} We define projections $ P_1 $, $ P_c $, and $ P_s $ on the modes concentrated around 
the wave vector $ k = e_1 $, on the modes concentrated around 
the wave vectors $ k = e_1 $ and $ k = - e_1 $, and on the linearly exponentially damped modes, respectively.
In detail we set 
$$ 
\widehat{P}_1(k) = \left\{ \begin{array}{cl} 1, & |k-e_1|\leq 1/10 , \\
0, & \textrm{else} ,\end{array} \right.
$$
$$ 
\widehat{P}_c(k) = \left\{ \begin{array}{cl} 1, & |k-e_1|\leq 1/10  \textrm{ or } |k+e_1|\leq 1/10, \\
0, & \textrm{else} ,\end{array} \right.
$$
and $ P_s = 1- P_c $. 
The projections $ P_1 $, $ P_c $, and $ P_s $ are bounded linear mappings from $ \mathcal{W}^r $ to $ \mathcal{W}^r $
for every $ r \geq 0 $.

\paragraph{Rescaling time in noise.}
Rescaling time in a Wiener process $ (\widehat{W}(k, t))_{t \geq 0}  $ for fixed $ k $ gives
the relation
$$
 \widehat{W}(k,t) =  \varepsilon^{-1} \widehat{W}(k,\varepsilon^2 t)
$$
and therefore 
$$ 
\widehat{\xi}(k,t) = \partial_t  \widehat{W}(k,t) =  \partial_t (\varepsilon^{-1} \widehat{W}(k, \varepsilon^2 t))
=  \varepsilon \widehat{\xi}(k, \varepsilon^2 t)
$$ 
in distributional sense.

\paragraph{Defining the critical part of the noise.}
In the following 
we assume that  $ (\alpha(k))_{k \in \mathcal{L}_{\varepsilon}} $ 
is of order $ \mathcal{O}(1) $  in $ L^1_{r_{SH}} $, i.e., we assume 
\begin{eqnarray*}
\|  \alpha\|_{L^1_{r_{SH}}} & = & \int_{\R^d}
| \alpha(k) | (1+|k|^2)^{r_{SH}/2}d\mu(k) \\
& = & \varepsilon^d\sum_{k \in \mathcal{L}_{\varepsilon}} | \alpha(k) | (1+|k|^2)^{r_{SH}/2}= C_{\alpha} = \mathcal{O}(1).
\end{eqnarray*}
From the above scaling relations we see that for fixed $ k \in \R $ we have to choose 
$$
 \widehat{W}_A(K,T) =  \varepsilon^{-1} \widehat{W}(k,\varepsilon^2 t),
$$
respectively 
$$
\widehat{\xi}_A(K,T) = \partial_T  \widehat{W}_A(K,T) =  \partial_t (\varepsilon^{-1} \widehat{W}(k, \varepsilon^2 t))
=  \varepsilon \widehat{\xi}(k, t)
$$
in distributional sense, where $ k = e_1 + \varepsilon K $ and $ T =  \varepsilon^2 t $.
We define $  \widehat{W}_A $  not for all $ k \in \mathcal{L}_{\varepsilon} $ in this way but only for those with
$$ 
k \in S_1 = \{ k \in \mathcal{L}_{\varepsilon} : | k-e_1  | < 1/10 \}.
$$
Therefore, to have  $ \zeta_A = \zeta_A(X,T) $, which is defined subsequently in \eqref{zetaA}, to be of order $ \mathcal{O}(1) $ in \eqref{GL} we have to set
$$
{\zeta}_1(x, t) = \varepsilon^d \sum_{k \in \mathcal{L}_{\varepsilon}\cap S_1} 
\varepsilon^2 \alpha(k) \widehat{\xi}(k,t) e^{ik \cdot x} ,
$$
where for given $ \alpha = \alpha(k) $ we set $ \alpha_A(K) = \alpha(e_1+ \varepsilon K) $ for $ | k-e_1 | < 1/10 $ and $ \alpha_A(K) = 0 $ elsewhere. 
Similarly, we set 
$$
{\zeta}_{-1}(x, t) = \varepsilon^d \sum_{k \in \mathcal{L}_{\varepsilon} \cap S_{-1}} 
\varepsilon^2 \alpha(k) \widehat{\xi}(k,t) e^{ik \cdot x} .
$$
Finally, we define the critical part of the noise $$ \zeta_c(x, t) =  {\zeta}_1(x, t) + {\zeta}_{-1}(x, t)  ,$$
 i.e. for $ |k-e_1| < 1/10 $ or $ |k+e_1| < 1/10 $ the critical part $ \zeta_c $ of the noise is of order 
$ \mathcal{O}(\varepsilon^2) $ w.r.t. to the $ \mathcal{W} $-norm in physical space.

\paragraph{Defining the stable  part of the noise.}
For the stable part, i.e., 
for 
$$ k \in S_s =  \{ k \in \mathcal{L}_{\varepsilon} :  |k-e_1| \geq 1/10 \textrm{ and } |k+e_1| \geq 1/10\} ,$$ 
the noise can  chosen 
to be larger, namely 
to be of order $ \mathcal{O}(\varepsilon^{\beta}) $ for a $ \beta > 1 $ w.r.t. to the $ \mathcal{W} $-norm in physical space.
Then we define the stable part $ \zeta_s $ of the noise as 
$$
{\zeta}_s(x, t) = \varepsilon^d \sum_{k \in \mathcal{L}_{\varepsilon} \cap S_s} 
\varepsilon^{\beta} \alpha(k) \widehat{\xi}(k,t) e^{ik \cdot x} .
$$

\paragraph{The spatial regularity.}
Since $ \alpha_A $ has a bounded support in Fourier space 
for fixed $ T $ the noise $ \zeta_A(T,\cdot) $ is arbitrarily smooth in space, i.e. 
$ \zeta_A(T,\cdot) \in \mathcal{W}_A^{r_A}$ for all
$ r_A  \geq 0 $. However, only w.r.t. to the $ \mathcal{W}_A $, i.e. for $  r_A = 0 $ the 
process $ \zeta_A $ is $ \mathcal{O}(1) $ and does not grow with $ \varepsilon
\to 0 $. This behavior is independent of the choice of $ r_{SH} \geq 0 $.
Therefore,  we choose $ r_{SH} = r_ A = 0 $ for the rest of the paper and come back 
to other choices in Section \ref{secdisc}.   Thus, we finally consider 
$$ 
{\zeta}(x, t) = {\zeta}_c(x, t) + {\zeta}_s(x, t) .
$$

\subsection{Derivation of the stochastic Ginzburg-Landau equation}

Then we set 
\begin{equation} \label{zetaA}
\widehat{\zeta}_A(K,T) = \varepsilon^{-2} P_c(e_1+ \varepsilon K)\widehat{\zeta}(e_1+ \varepsilon K,\varepsilon^2 t).
\end{equation}
The stochastic Ginzburg-Landau equation is then given by the deterministic Ginzburg-Landau 
equation \eqref{GLdet} plus the noise term $ {\zeta}_A $.
Therefore, if we define 
$$
\textrm{Res}_{stoch}(u) = - \partial_t u - (1+ \partial_{x_1}^2)^2 u + 4 \Delta_{x_{\perp}} u + \varepsilon^2 u - u^3 +\zeta
$$
we obtain 
$$
\textrm{Res}_{stoch}(u) = \textrm{Res}_{det}(u)  + P_s \zeta.
$$

\subsection{The approximation result}

Our approximation result is then as follows 
\begin{theorem} \label{mainth}
For a fixed but arbitrary small $ \delta' > 0 $,
for all $ \beta > 1 $, $ \delta > 0 $, and $ C_1 > 0 $  there exist $ \varepsilon_0 > 0 $, $ C_2 > 0 $  such that for all $ \varepsilon \in (0,\varepsilon_0) $ the following holds. Let 
$ B $ 
be a solution 
of \eqref{vocnewGL} and 
$ Z_A $ be defined in \eqref{Zeqvo}
with 
\begin{equation} \label{ass3}
\sup_{\tau \in [0,T_0]}\|Z_A(\cdot,\tau) \|_{\mathcal{W}_A^{1}} + \sup_{\tau \in [0,T_0]}\|{B}(\cdot,\tau) \|_{\mathcal{W}_A^{3-\delta'}} \leq C_1 .
\end{equation}
Then there are solutions 
$ u $ 
of  \eqref{SH} with 
$$ \mathcal{P}(\sup_{t \in [0,T_0/\varepsilon^2]} \sup_{x \in \R}|u(x,t)-(\varepsilon (Z_A(\varepsilon x, \varepsilon^2 t)+ B(\varepsilon x, \varepsilon^2 t)) e^{i x_1}+ c.c) |  
\leq C_2 \varepsilon^\beta)   >   1 - \delta.
$$
\end{theorem}
\begin{remark}{\rm 
Assumption \eqref{ass3} can obviously be  satisfied in the sense of Corollary \ref{cor34a} and Corollary \ref{cor34}.
}
\end{remark}
\begin{remark}{\rm 
Theorem \ref{mainth} is an improvement of the existing results in the literature, cf. \cite{BHP05}, even in the 1D case, 
in the sense that the noise in the stable part can be $ \mathcal{O}(\varepsilon^{\beta}) $ with $ \beta > 1 $ instead 
of $ \mathcal{O}(\varepsilon^{2}) $.
}
\end{remark}

\subsection{The equations for the error} 

We write an exact solution $ u $ of \eqref{SH} as sum of the GL approximation 
$ \varepsilon \Psi $ and an error $  \varepsilon^{\beta} R $, i.e., $ u = \varepsilon \Psi + \varepsilon^{\beta} R $ for $ \beta \in (1,2] $
and find that $ R $ satisfies 
\begin{eqnarray*}
\partial_t R & = & - (1+\partial_{x_1}^2)^2 R + 4 \Delta_{x_{\perp}} R + \varepsilon^2 R 
\\ && -
3 \varepsilon^2 \Psi^2 R 
- 3 \varepsilon^{1+\beta} \Psi R^2 - \varepsilon^{2 \beta} R^3 + \varepsilon^{-\beta} \textrm{Res}_{det}(\varepsilon \Psi)+ \varepsilon^{-\beta} P_s \zeta .
\end{eqnarray*}
We split $ R = R_c + R_s $ where $ R_c = P_c R $ and $ R_s = P_s R $. This gives the system
\begin{eqnarray*}
\partial_t R_c & = & - (1+\partial_{x_1}^2)^2 R_c + 4 \Delta_{x_{\perp}} R_c + \varepsilon^2 R_c 
\\ && + P_c(-
3 \varepsilon^2 \Psi^2 R 
- 3 \varepsilon^{1+\beta} \Psi R^2 - \varepsilon^{2 \beta} R^3  ) + \varepsilon^{-\beta} P_c \textrm{Res}_{det}(\varepsilon \Psi) , 
\\
\partial_t R_s & = & - (1+\partial_{x_1}^2)^2 R_s + 4 \Delta_{x_{\perp}} R_s + \varepsilon^2 R_s 
\\ && + P_s(-
3 \varepsilon^2 \Psi^2 R 
- 3 \varepsilon^{1+\beta} \Psi R^2 - \varepsilon^{2 \beta} R^3 ) + \varepsilon^{-\beta} P_s \textrm{Res}_{det}(\varepsilon \Psi)+ \varepsilon^{-\beta} P_s \zeta .
\end{eqnarray*}
The final step is the splitting of $ R_s = R_B + R_Z $ into a regular part $ R_B $ and into 
an Ornstein-Uhlenbeck process $ R_Z $. Therefore, we finally consider 
\begin{eqnarray*}
\partial_t R_c & = & - (1+\partial_{x_1}^2)^2 R_c + 4 \Delta_{x_{\perp}} R_c + \varepsilon^2 R_c 
\\ && + P_c(-
3 \varepsilon^2 \Psi^2 R 
- 3 \varepsilon^{1+\beta} \Psi R^2 - \varepsilon^{2 \beta} R^3 ) + \varepsilon^{-\beta} P_c \textrm{Res}_{det}(\varepsilon \Psi) , 
\\
\partial_t R_B & = & - (1+\partial_{x_1}^2)^2 R_B + 4 \Delta_{x_{\perp}} R_B + \varepsilon^2 R_B 
+ \varepsilon^2 R_Z
\\ && + P_s(-
3 \varepsilon^2 \Psi^2 R 
- 3 \varepsilon^{1+\beta} \Psi R^2 - \varepsilon^{2 \beta} R^3  ) + \varepsilon^{-\beta} P_s \textrm{Res}_{det}(\varepsilon \Psi), \\
\partial_t R_Z & = & - (1+\partial_{x_1}^2)^2 R_Z + 4 \Delta_{x_{\perp}} R_Z + \varepsilon^{-\beta} P_s \zeta .
\end{eqnarray*}

\subsection{The Ornstein-Uhlenbeck process $ R_Z $}
\label{sec55}

We start with the third equation 
$$ 
\partial_t R_Z  =  - (1+\partial_{x_1}^2)^2 R_Z + 4 \Delta_{x_{\perp}} R_Z + \varepsilon^{-\beta} P_s \zeta 
$$
which again will be solved with the variation of constant formula.
Using $  \widehat{R}_Z(k,0) = 0 $ we find for $ k \in \mathcal{L}_{\varepsilon}  $ that
\begin{eqnarray*}
\widehat{R}_Z(k,t) & = & \varepsilon^{-\beta} \int_0^t  e^{\lambda(k) (t-\tau)}  \widehat{P}_s(k) \alpha(k) \widehat{\xi}(k,\tau) d\tau
\\ & = & 
\varepsilon^{-\beta}  \int_0^t  e^{\lambda(k) (t-\tau)} \widehat{P}_s(k) \alpha(K) d\widehat{W}(k,\tau)  ,
\end{eqnarray*}
where $ \lambda(k)  = -(1- |k_1|^2)^2 - 4 |k_{\perp}|^2 $ and  the integral is a stochastic integral in the It\^{o}-sense. Following the calculations of Section \ref{sec3} 
and using the fact that $ \lambda(k) $ is strictly negative for $ R_s $
such that 
$$
\int_0^t e^{\lambda(k)(t-\tau)} \widehat{P}_s(k) d\tau \leq C/\lambda(k) 
$$
uniformly for $ t \geq 0 $
we obtain 
$$
\text{P}[\sup_{\tau \in [0,t]}\|\widehat{R}_Z(\cdot,\tau) \|_{L^1_r} \geq c] \leq 
 \frac{C}{c^2}\sup_{k \in \R^d} \frac{(1+|k|^2)^{r-r_{SH}}}{\lambda(k) }
\int_{\R^d} \alpha(k) (1+|k|^2)^{r_{SH}/2}  d\mu(k).
$$
Hence, we finally find 
\begin{lemma} \label{lemshRZ}
Suppose $ r \leq  r_{sh} +1 $.
Then for all $ \delta_0 > 0 $ there 
exists a $ c > 0 $ such that 
$$
\text{P}[\sup_{\tau \in [0,t]}\|\widehat{R}_Z(\cdot,\tau) \|_{L^1_r} \geq c] \leq \delta_0
$$
for all $ t \geq 0 $.
\end{lemma}
In the following we only need this estimate for $ r = 0 $.
\begin{corollary} \label{coroRZ}
Let $ r_{sh} = 0  $.
Then for all $ \delta_0 > 0 $ there 
exists a $ c > 0 $ such that 
$$
\text{P}[\sup_{\tau \in [0,t]}\|{R}_Z(\cdot,\tau) \|_{\mathcal{W}} \geq c] \leq \delta_0
$$
for all $ t \geq 0 $.
\end{corollary}

\subsection{Estimates for the residual}

\label{res56}

In this section we provide estimates for the residual
\begin{eqnarray*}
\textrm{Res}_{det}(u_A) & = &  4 i \varepsilon^4 \partial_{X_1}^3  A_1 e^{ix_1} - 4 i \varepsilon^4 \partial_{X_1}^3 A_{-1} e^{-ix_1}  \\ && 
-  \varepsilon^5 \partial_{X_1}^4  A_1 e^{ix_1} -  \varepsilon^5 \partial_{X_1}^4 A_{-1} e^{-ix_1} 
\\ && - \varepsilon^3 A_1^3 e^{3ix_1}-  \varepsilon^3 A_{-1}^3 e^{-3ix_1}  .
\end{eqnarray*}	
The equations for the error will be solved with the help of the variation of constant formula. Therefore, we also need the  estimate 
$$ 
\textrm{RES}_{j}(t) =  \int_0^t e^{(- (1+\partial_{x_1}^2)^2  + 4 \Delta_{x_{\perp}})(t - \tau)} P_j (\textrm{Res}_{det}(\varepsilon \Psi))(\tau)
d \tau 
$$
for $ j = c,s $ in $ \mathcal{W} $.
We have 
\begin{lemma} \label{lem54}
For all $ \delta, \delta_0 > 0 $ 
there exist $ C,\varepsilon_0 > 0 $ such that for all $ \varepsilon \in (0,\varepsilon_0) $ we have
\begin{eqnarray*}
P(\|\textrm{RES}_{c}(\cdot,t) \|_{\mathcal{W}} \leq  C \varepsilon^{2-\delta}) + P(\|\textrm{RES}_{s}(\cdot,t) \|_{\mathcal{W}} \leq  C \varepsilon^{3-\delta})& \geq & 1 - \delta_0.
\end{eqnarray*}	
\end{lemma}
\noindent
{\bf Proof.}  
The estimates for the deterministic part of the stochastic 
residual are exactly the same as the estimates in the pure deterministic case.
The ingredients which we need for estimating this part   are as follows.

i) With probability $ 1- \delta_0 $ we have that $ A_1 $ is uniformly $ \mathcal{O}(1) $-bounded 
in $ \mathcal{W}_{r_A} $.

iia) For such $ A_1 $ we obviously have 
$$ 
\|P_s(- \varepsilon^3 A_1^3 e^{3ix_1}-  \varepsilon^3 A_{-1}^3 e^{-3ix_1}) \|_{\mathcal{W}} \leq C \varepsilon^3.
$$ 

iiia) For estimating 
$$ 
P_c (- \varepsilon^3 A_1^3 e^{3ix_1}-  \varepsilon^3 A_{-1}^3 e^{-3ix_1}) 
$$
we use that $ A_1^3 \in \mathcal{W}_{r_A} $ and that $ \varepsilon^3 A_1^3 e^{3ix_1} $
is concentrated in Fourier space at $ 3 e_1 $. The exact estimate is given below.

iv) The term
$$
( 4 i \varepsilon^4 \partial_{X_1}^3   -  \varepsilon^5 \partial_{X_1}^4)  A_1 e^{ix_1} 
$$
is given in Fourier space by 
$$
 (- 4  \varepsilon^4 K_1^3 - \varepsilon^5 K_1^4) \widehat{A}_1(K), \qquad \textrm{with} \qquad  k = e_1+ \varepsilon K 
= e_1+ \varepsilon (K_1,K_{\perp}) .
$$ 
The pre-factor can be interpreted as 
\begin{equation} \label{star} 
- 4  \varepsilon^4 K_1^3 - \varepsilon^5 K_1^4 = \varepsilon(\lambda(k) + 4 (k_1-1)^2+ 4 |k_{\perp}|^2).
\end{equation}
To estimate this term we use that $ A_1 = B + Z_A $ with 
$ B $ 
and
$ Z_A $ 
are bounded according to Corollary \ref{cor34} and Corollary \ref{cor34a}.

To estimate
$$ 
P_1 (4 i \varepsilon^4 \partial_{X_1}^3   -  \varepsilon^5 \partial_{X_1}^4)  B e^{ix_1} 
$$
we use that 
\begin{equation} \label{why} 
|\widehat{P}_1(k)(\lambda(k) + 4 (k_1-1)^2+ 4 |k_{\perp}|^2)) | \leq C |k_1-1|^{3}
\end{equation}
due to the compact support of $ \widehat{P}_1$.
For the same reason we can estimate
\begin{equation} \label{why1} 
|k_1-1|^{3} \leq 
C |k_1-1|^{3-\delta}
\end{equation}
which we have to use at one point due to the low regularity of 
$ B $.
Hence, for all $ t \in [0,T_0/\varepsilon^2] $ we  get 
with \eqref{star} and \eqref{why} that
\begin{eqnarray*}
\lefteqn{ \| \int_0^t e^{(- (1+\partial_{x_1}^2)^2  + 4 \Delta_{x_{\perp}})(t - \tau)} P_1 (( 4 i \varepsilon^4 \partial_{X_1}^3   -  \varepsilon^5 \partial_{X_1}^4)  B(\varepsilon \cdot, \varepsilon^2 \tau) e^{ix_1} ) d\tau \|_{\mathcal{W}}} \\
& =&\int_{R^d} | \int_0^t e^{\lambda(k)(t-\tau)}(\lambda(k) + 4 (k_1-1)^2)+ 4 |k_{\perp}|^2) \\ && \qquad \qquad \times\varepsilon
\varepsilon^{-d} \widehat{P}_1(k) \widehat{B}
(\frac{k-e_1}{\varepsilon}, \varepsilon^2 \tau) d\tau| d\mu(k) \\
& = & \int_{R^d} | \int_0^t e^{\lambda(e_1+ \varepsilon K)(t-\tau)} (\lambda(e_1+ \varepsilon K) + 4 \varepsilon^2 K_1^2 + 4 |k_{\perp}|^2)) \\ && \qquad \qquad \times\varepsilon \widehat{P}_1(e_1+ \varepsilon K)
 \widehat{B}
(K,\varepsilon^2 \tau) d\tau | d\mu_A(K) \\
& \leq & C \int_{R^d} |\int_0^t e^{\lambda(e_1+ \varepsilon K)(t-\tau)} \varepsilon^3 K^3 \varepsilon
\widehat{P}_1(e_1+ \varepsilon K) \widehat{B}
(K,\varepsilon^2 \tau) d\tau| d\mu_A(K) 
\\
& \leq & C \int_{R^d} |\int_0^t e^{\lambda(e_1+ \varepsilon K)(t-\tau)} \varepsilon^{\delta} K^{\delta} \varepsilon \varepsilon^{3-\delta} K^{3-\delta} \widehat{P}_1(e_1+ \varepsilon K)
 \widehat{B}
(K,\varepsilon^2 \tau) d\tau | d\mu_A(K)\\
& \leq & C  \int_0^t \sup_{K \in \R^d}| e^{-\varepsilon^2 K^2 (t-\tau)} \varepsilon^{\delta} K^{\delta} | d\tau   \cdot \varepsilon^{4-\delta} 
 \sup_{T \in [0,T_0]}\| \widehat{B}
(\cdot,T) \|_{L^1_{3-\delta}}   \\
& \leq & C  \int_0^t  (t-\tau)^{-\delta/2} d\tau   \varepsilon^{4-\delta} 
 \sup_{T \in [0,T_0]} \| \widehat{B}
(\cdot,T) \|_{L^1_{3-\delta}}   \\
& \leq & C  | (T_0/\varepsilon^2-\tau)^{1-\delta/2} |_{\tau=0}^{T_0/\varepsilon^2}  | \cdot \varepsilon^{4-\delta} 
 \sup_{T \in [0,T_0]} \| \widehat{B}
(\cdot,T) \|_{L^1_{3-\delta}}   \\
& \leq & C \varepsilon^2
\sup_{T \in [0,T_0]} \| \widehat{B}
(\cdot,T) \|_{L^1_{3-\delta}} .
\end{eqnarray*}
Similarly, for all $ t \in [0,T_0/\varepsilon^2] $ we  get 
with \eqref{star}, \eqref{why}, and \eqref{why1} that
\begin{eqnarray*}
\lefteqn{ \| \int_0^t e^{(- (1+\partial_{x_1}^2)^2  + 4 \Delta_{x_{\perp}})(t - \tau)} P_1 (( 4 i \varepsilon^4 \partial_{X_1}^3   -  \varepsilon^5 \partial_{X_1}^4)  Z_A(\varepsilon \cdot, \varepsilon^2 \tau) e^{ix_1} ) d\tau \|_{\mathcal{W}}} \\
& =&\int_{R^d} | \int_0^t e^{\lambda(k)(t-\tau)}(\lambda(k) + 4 (k_1-1)^2+ 4|k_{\perp}|^2))
\\ && \qquad \qquad \times \varepsilon
\varepsilon^{-d} \widehat{P}_1(k) \widehat{Z}_A
(\frac{k-e_1}{\varepsilon}, \varepsilon^2 \tau) d\tau| d\mu(k) \\
& = & \int_{R^d} | \int_0^t e^{\lambda(e_1+ \varepsilon K)(t-\tau)} (\lambda(e_1+ \varepsilon K) +4 \varepsilon^2 K_1^2+ 4|k_{\perp}|^2)) 
\\ && \qquad \qquad \times \varepsilon  \widehat{P}_1(e_1+ \varepsilon K)
 \widehat{Z}_A
(K,\varepsilon^2 \tau) d\tau | d\mu_A(K) \\
& \leq & C \int_{R^d} |\int_0^t e^{\lambda(e_1+ \varepsilon K)(t-\tau)} \varepsilon^{3-\delta} K^{3-\delta} 
\\ && \qquad \qquad \times \varepsilon \widehat{P}_1(e_1+ \varepsilon K) \widehat{Z}_A
(K,\varepsilon^2 \tau) d\tau| d\mu_A(K) 
\\
& \leq & C \int_{R^d} |\int_0^t e^{\lambda(e_1+ \varepsilon K)(t-\tau)} \varepsilon^{2-\delta/2} K^{2-\delta/2} \varepsilon \varepsilon^{1-\delta/2} K^{1-\delta/2} \\ && \qquad \qquad \times \widehat{P}_1(e_1+ \varepsilon K)
 \widehat{Z}_A
(K,\varepsilon^2 \tau) d\tau | d\mu_A(K)\\
& \leq & C  \int_0^t \sup_{K \in \R^d}| e^{-\varepsilon^2 K^2 (t-\tau)} \varepsilon^{2-\delta/2} K^{2-\delta/2} | d\tau   \cdot \varepsilon^{2-\delta/2} 
 \sup_{T \in [0,T_0]}\| \widehat{Z}_A
(\cdot,T) \|_{L^1_{1-\delta/2}}   \\
& \leq & C  \int_0^t  (t-\tau)^{-(1-\delta/2)} d\tau   \cdot \varepsilon^{2-\delta/2} 
 \sup_{T \in [0,T_0]} \| \widehat{Z}_A
(\cdot,T) \|_{L^1_{1-\delta/2}}   \\
& \leq & C  | (T_0/\varepsilon^2-\tau)^{\delta/2} |_{\tau=0}^{T_0/\varepsilon^2}  | \cdot \varepsilon^{2-\delta/2} 
 \sup_{T \in [0,T_0]} \| \widehat{Z}_A
(\cdot,T) \|_{L^1_{1-\delta/2}}   \\
& \leq & C \varepsilon^{2-3\delta/2} 
\sup_{T \in [0,T_0]} \| \widehat{Z}_A
(\cdot,T) \|_{L^1_{1-\delta/2}} .
\end{eqnarray*}

iib) Due to the support of $ \widehat{P}_s $ the associated part of the linear semigroup is damped with a rate $ e^{- \sigma (t - \tau)} $ for a $ \sigma > 0 $. Thus, we obtain
\begin{eqnarray*}
\lefteqn{ \| \int_0^t e^{(- (1+\partial_{x_1}^2)^2  + 4 \Delta_{x_{\perp}})(t - \tau)} P_s (- \varepsilon^3 A_1(\varepsilon \cdot, \varepsilon^2 \tau)^3 e^{3ix_1}-  \varepsilon^3 A_{-1}(\varepsilon \cdot, \varepsilon^2 \tau)^3 e^{-3ix_1}) d\tau \|_{\mathcal{W}} }\\
& \leq & C \int_0^t e^{- \sigma (t - \tau)} d\tau \sup_{\tau \in [0,t]} 
\|P_s(- \varepsilon^3 A_1(\varepsilon \cdot, \varepsilon^2 \tau)^3 e^{3ix_1}-  \varepsilon^3 A_{-1}(\varepsilon \cdot, \varepsilon^2 \tau)^3 e^{-3ix_1}) \|_{\mathcal{W}} 
\\ & \leq & C \varepsilon^3.
\end{eqnarray*}

iiib) To estimate
$$\| \int_0^t e^{(- (1+\partial_{x_1}^2)^2  + 4 \Delta_{x_{\perp}})(t - \tau)} P_c (- \varepsilon^3 A_1(\varepsilon \cdot, \varepsilon^2 \tau)^3 e^{3ix_1}-  \varepsilon^3 A_{-1}(\varepsilon \cdot, \varepsilon^2 \tau)^3 e^{-3ix_1})  d\tau \|_{\mathcal{W}} 
$$
we use that $ P_c = P_1 + P_{-1} $ and that it is sufficient to estimate the $ P_1 $-part.
we find 
\begin{eqnarray*}
\lefteqn{ \| \int_0^t e^{(- (1+\partial_{x_1}^2)^2  + 4 \Delta_{x_{\perp}})(t - \tau)} P_1 ( \varepsilon^3 A_1(\varepsilon \cdot, \varepsilon^2 \tau)^3 e^{3ix_1}+  \varepsilon^3 A_{-1}(\varepsilon \cdot, \varepsilon^2 \tau)^3 e^{-3ix_1})  d\tau \|_{\mathcal{W}} }
\\ & = & \int_{\R^d}|
 \int_0^t e^{\lambda(k)(t - \tau)} \\ && \qquad  \times  \widehat{P}_1(k) ( \varepsilon^3 \varepsilon^{-d} \widehat{A}_1^{*3} (\frac{k-3e_1}{\varepsilon}, \varepsilon^2 \tau)+  \varepsilon^3 \varepsilon^{-d} \widehat{A}_{-1}^{*3} (\frac{k+3e_1}{\varepsilon}, \varepsilon^2 \tau))  d\tau
|d\mu(k) \\
\\ & \leq  &
 \int_0^t \sup_{k \in \R^d} |e^{\lambda(k)(t - \tau)}| \\ && \qquad \times  \int_{\R^d}| \widehat{P}_1(k) ( \varepsilon^3 \varepsilon^{-d} \widehat{A}_1^{*3} (\frac{k-3e_1}{\varepsilon}, \varepsilon^2 \tau)+  \varepsilon^3 \varepsilon^{-d} \widehat{A}_{-1}^{*3} (\frac{k+3e_1}{\varepsilon}, \varepsilon^2 \tau)) 
|d\mu(k) d\tau
\\ & \leq  &
 \int_0^t e^{\varepsilon^2 (t - \tau)}  (\int_{\R^d}| \widehat{P}_1(k+3 e_1) ( \varepsilon^3 \varepsilon^{-d} \widehat{A}_1^{*3} (\frac{k}{\varepsilon}, \varepsilon^2 \tau)
|d\mu(k) 
\\ && \qquad \qquad \qquad+ 
 \int_{\R^d}| \widehat{P}_1(k-3e_1)   \varepsilon^3 \varepsilon^{-d} \widehat{A}_{-1}^{*3} (\frac{k}{\varepsilon}, \varepsilon^2 \tau)) 
|d\mu(k)) d\tau
\\ & \leq  & T_0 e^{T_0} \varepsilon^{-2}
\sup_{\tau \in [0,T_0/\varepsilon^2]} \int_{\R^d}| \widehat{P}_1(\varepsilon K +3 e_1) ( \varepsilon^3 \widehat{A}_1^{*3} (K, \varepsilon^2 \tau)
|d\mu_A(K) 
\\ && \qquad + T_0 e^{T_0} \varepsilon^{-2}
\sup_{\tau \in [0,T_0/\varepsilon^2]} 
 \int_{\R^d}| \widehat{P}_1(\varepsilon K-3e_1)   \varepsilon^3 \widehat{A}_{-1}^{*3} (K, \varepsilon^2 \tau)) 
|d\mu_A(K)
\\ & \leq  & T_0 e^{T_0} \varepsilon^{-2} \sup_{K \in \R^d} | \widehat{P}_1(\varepsilon K +3 e_1) (1+K^2)^{-1/2+\delta/2} |
\\ && \qquad \times  \sup_{\tau \in [0,T_0/\varepsilon^2]}
 \int_{\R^d}|  ( \varepsilon^3 \widehat{A}_1^{*3} (K, \varepsilon^2 \tau)  |(1+K^2)^{1/2-\delta/2}  
d\mu_A(K) 
\\ &&  +  T_0 e^{T_0} \varepsilon^{-2}
 \sup_{K \in \R^d} |  \widehat{P}_1(\varepsilon K-3e_1) (1+K^2)^{-1/2+\delta/2} |
\\ && \qquad \times 
\sup_{\tau \in [0,T_0/\varepsilon^2]} \int_{\R^d}|   \varepsilon^3 \widehat{A}_{-1}^{*3} (K, \varepsilon^2 \tau)) |(1+K^2)^{1/2-\delta/2}  
|d\mu_A(K)
\\ & \leq & C T_0 e^{T_0}  \varepsilon^{2-\delta} \sup_{\tau \in [0,T_0]} \| A_1(\cdot,T) \|_{\mathcal{W}_{A}^{1-\delta}}^3,
\end{eqnarray*}
since $  \sup_{K \in \R^d} |  \widehat{P}_1(\varepsilon K\pm 3e_1) (1+K^2)^{-1/2+\delta/2} | = \mathcal{O}(\varepsilon^{1-\delta}) $.

v) We have $ P_s (Z_A e^{i x_1}) = 0 $ such that the last estimate is given by
\begin{eqnarray*}
&& \| \int_0^t e^{(- (1+\partial_{x_1}^2)^2  + 4 \Delta_{x_{\perp}})(t - \tau)} P_s (( 4 i \varepsilon^4 \partial_{X_1}^3   -  \varepsilon^5 \partial_{X_1}^4)  B(\cdot,\varepsilon^2 \tau) e^{ix_1} ) d\tau \|_{\mathcal{W}} \\
&  = & \int_{R^d} | \int_0^t e^{\lambda(k)(t - \tau)} \widehat{P}_s(k) ((-4 \varepsilon^4 K_1^3   -  \varepsilon^5 K_1^4)\\ && \qquad \qquad \times  \varepsilon^{-d} \widehat{B}(\frac{k-e_1}{\varepsilon},\varepsilon^2 \tau)) d\tau |d\mu(k) \\
&  = & \int_{R^d} | \int_0^t e^{\lambda(k)(t - \tau)} \widehat{P}_s(k) ((-4 \varepsilon (k_1-1)^3   -  \varepsilon (k_1-1)^4) \\ && \qquad \qquad \times \varepsilon^{-d} \widehat{B}(\frac{k-e_1}{\varepsilon},\varepsilon^2 \tau)) d\tau |d\mu(k) \\
&  \leq &   \int_0^t \sup_{k \in \R^d} |e^{\lambda(k)(t - \tau)} \widehat{P}_s(k) ((4 \varepsilon |k_1-1|   +  \varepsilon |k_1-1|^2) \\ && \qquad \qquad \times 
 \int_{R^d} |\varepsilon^{-d} \widehat{B}(\frac{k-e_1}{\varepsilon},\varepsilon^2 \tau)) |k_1-1|^2   |d\mu(k)  d\tau\\
 &  \leq & \varepsilon \int_0^t C e^{-\sigma(t-\tau)}  (t-\tau)^{-1/2} d\tau
\\ && \qquad \qquad \times  \sup_{\tau \in [0,T_0/\varepsilon^2]} \int_{R^d} |\varepsilon^{-d} \widehat{B}(\frac{k-e_1}{\varepsilon},\varepsilon^2 \tau)) |k_1-1|^2   |d\mu(k)  \\
 &  \leq & C \varepsilon  \sup_{\tau \in [0,T_0/\varepsilon^2]} \int_{R^d} |\widehat{B}(K,\varepsilon^2 \tau)) \varepsilon^2 |K|^2   |d\mu_A(K) 
\\ 
 &  \leq & C \varepsilon^3 \sup_{T \in [0,T_0]}\|B(\cdot,T) \|_{\mathcal{W}_{A}^2}.
\end{eqnarray*}
\qed

%
%
%
%
%
%
%
%
%
%
%
%
%

\subsection{The error estimates}

It remains  to estimate  $ R_c $ and $ R_s $ as solutions of the first two equations
\begin{eqnarray*}
\partial_t R_c & = & - (1+\partial_{x_1}^2)^2 R_c + 4 \Delta_{x_{\perp}} R_c + \varepsilon^2 R_c 
\\ && + P_c(-
3 \varepsilon^2 \Psi^2 R 
- 3 \varepsilon^{1+\beta} \Psi R^2 - \varepsilon^{2 \beta} R^3  ) + \varepsilon^{-\beta} P_c \textrm{Res}(\varepsilon \Psi) , 
\\
\partial_t R_B & = & - (1+\partial_{x_1}^2)^2 R_B + 4 \Delta_{x_{\perp}} R_B + \varepsilon^2 R_B 
+ \varepsilon^2 R_Z
\\ && + P_s(-
3 \varepsilon^2 \Psi^2 R 
- 3 \varepsilon^{1+\beta} \Psi R^2 - \varepsilon^{2 \beta} R^3  ) + \varepsilon^{-\beta} P_s \textrm{Res}(\varepsilon \Psi), 
\end{eqnarray*}
with $ R = R_c + R_B + R_Z $ and $ \Psi =  \varepsilon (B+Z_A) e^{ix_1} + c.c. $.
We recall that we have chosen
$ r_{SH} = 0 $ and $ r_A = 0 $. 
Therefore, we know that $ \Psi $ satisfies 
$$ 
\sup_{t \in [0,T_0/\varepsilon^2]} \| \Psi(\cdot,t) \|_{\mathcal{W}^{1}} =  \mathcal{O}(1) .
$$ 
Moreover, we have that $ R_Z $ 
satisfies 
$$ 
\sup_{t \in [0,T_0/\varepsilon^2]} \| R_Z(\cdot,t) \|_{\mathcal{W}} =  \mathcal{O}(1) 
$$ 
due to 
Corollary \ref{coroRZ}.
The system for $ R_c $ and $ R_B $ is solved with the help of the variation of constant formula
\begin{eqnarray*}
&& R_c(t)  =  \int_0^t e^{(- (1+\partial_{x_1}^2)^2 + 4 \Delta_{x_{\perp}})(t -\tau)}
\\ && \qquad \times (\varepsilon^2 R_c 
+ P_c(-
3 \varepsilon^2 \Psi^2 R 
- 3 \varepsilon^{1+\beta} \Psi R^2 - \varepsilon^{2 \beta} R^3  )(\tau) d \tau + \varepsilon^{-\beta}  \textrm{RES}_c(t) , 
\\
&&R_B(t)  =  \int_0^t e^{(- (1+\partial_{x_1}^2)^2 + 4 \Delta_{x_{\perp}})(t -\tau)}
\\ && \qquad \times ( \varepsilon^2 R_B 
+ \varepsilon^2 R_Z
 + P_s(-
3 \varepsilon^2 \Psi^2 R 
- 3 \varepsilon^{1+\beta} \Psi R^2 - \varepsilon^{2 \beta} R^3  ) + \varepsilon^{-\beta} \textrm{RES}_s(t).
\end{eqnarray*}
Lemma \ref{lem54} implies the uniform $ \mathcal{O}(1) $-boundedness of the residual terms 
$ \varepsilon^{-\beta} \textrm{RES}_c(t) $ and $  \varepsilon^{-\beta} \textrm{RES}_s(t) $ 
in $ \mathcal{W} $.

%

Like in the deterministic case we introduce 
$$ 
S(t) = \sup_{0 \leq \tau \leq t} \| R_c(t) \|_{\mathcal{W}} + \sup_{0 \leq \tau \leq t} \| R_B(t) \|_{\mathcal{W}}.
$$
For all $ t \in [0,T_0/\varepsilon^2] $ we obtain the inequality 
$$ 
S(t) \leq \int_0^t C_1 \varepsilon^2 S(\tau) + C_2 \varepsilon^3 S^2(\tau) + C_3 \varepsilon^4 S^3(\tau)
d \tau + C_{res},
$$ 
with constants $ C_1,C_2,C_3 $ and $ C_{res} $ independent of $ 0 < \varepsilon^2 \ll 1 $, the quantity $ S(t) $, and $ t \geq 0 $. 
Therefore, there exists an $ \varepsilon > 0 $ such that for all $ \varepsilon \in (0,\varepsilon_0) $ a simple application of Gronwall's inequality to the variation of constant formula yields that 
$ R_c $ and 
$ R_B $ are uniform $ \mathcal{O}(1) $-bounded in $ \mathcal{W} $.
This proves Theorem \ref{mainth}. \qed

\section{Outlook and discussion}
\label{secdisc}

In this section, we make remarks about  the handling of quadratic nonlinearities, the handling of  noise of smaller 
magnitude, and the handling of space-fractional  Swift-Hohenberg models.

\begin{remark}{\rm
The handling of quadratic nonlinearities is more involved. A possible  toy problem is 
\begin{equation} \label{quad}
\partial_t u = - (1+\partial_{x_1}^2)^2 u + 4 \Delta_{x_{\perp}} u + \alpha u -u \nabla \cdot u + \zeta 	.
\end{equation}
The ansatz for the derivation of the Ginzburg-Landau equation is 
\begin{equation}\label{ansatz1}
u(x,t) = \varepsilon A_1(X,T) e^{ix_1} + \varepsilon^2 A_2(X,T) e^{2ix_1} + \varepsilon^2 A_0(X,T) 
+ c.c. .
\end{equation}
Although in the deterministic case the separation of the critical and stable modes plays a fundamental role in 
proving the approximation result, cf. \cite{vH91,Schn94a,Schn94b},  the scaling of the error parts $ R_c $ and $ R_B $ in the cubic and 
quadratic case would differ. Moreover, in the residual the terms $ \partial_T A_2 $ and $ \partial_T A_0 $
must be estimated, which leads to stochastic integrals in the sense of It\^{o}. 
Due to the different scalings and the more involved residual estimates  the 
 handling 
of quadratic terms will be the subject of a forthcoming paper.}
\end{remark}

%
%
%
%

\begin{remark}{\rm
An obvious remark is that if $  \zeta_c(\cdot,t)  = \mathcal{O}(\varepsilon^{1+\beta}) $ with $ \beta > 1 $ 
in the above sense, then the deterministic  Ginzburg-Landau  equation 
$$
\partial_T A = 4 \Delta_X A  + A - 3 A |A|^2 
$$ 
occurs.
Since the solutions of this Ginzburg-Landau equation are arbitrarily smooth,
the regularity problems do not play a role in this case and the proof of the approximation result 
is only a slight modification of the deterministic case
if the noise is so small that it does not appear in the Ginzburg-Landau equation. 
}
\end{remark}

%
%
%
%
%
%

\begin{remark}{\rm
Space-fractional equations recently attracted a lot of interest.
A Ginzburg-Landau approximation result for the space-fractional Swift-Hohenberg equation 
can be found in  \cite{KT25}. 
For the derivation of the approximation equations the space-fractional character 
only plays a role  at the wave number $ k = 0 $.
Therefore, 
since  $ k = 0 $ does not appear in the derivation of the Ginzburg-Landau  equation,
the space-fractional character plays no role 
for the Swift-Hohenberg equation and so our approximation result also applies for
$$
\partial_t u = - (1-(-\partial_{x_1}^{2})^{\theta/2})^2 u + 4 \Delta_{x_{\perp}} u + \alpha u -u^3 + \zeta 
$$ 
for $ \vartheta > 0 $ sufficiently large. There are only two points, one has to be careful about. Firstly, for 
the validity of Lemma \ref{lemshRZ}, respectively Corollary \ref{coroRZ}, one needs $ \theta > 1/4 $ since now 
$ \lambda(k) = -(1-k_1^2)^{2\theta} - 4 |k_{\perp}|^2 $ which decays as $ - |k_1|^{4 \theta} - 4 |k_{\perp}|^2 $ for $ |k| \to \infty $.
Secondly, for the estimates  in Subsection \ref{res56}, 
Taylor expansions of $ \lambda $ at $ (k_1,k_{\perp}) = (\pm 1,0)$ have to be used. However, this will not lead to any 
further restriction on $ \vartheta $.} 
\end{remark}

\bibliographystyle{alpha}
\bibliography{GLbib}

\end{document}